\documentclass[12pt]{amsart}

\usepackage[english]{babel}
\usepackage{vmargin,graphicx,amsmath,amssymb,color,subfigure,enumerate,csquotes, dsfont, mathrsfs,color}
\numberwithin{equation}{section}

\usepackage{vmargin,color}
\definecolor{gr}{rgb}   {0.,   0.69,   0.23 }
\definecolor{bl}{rgb}   {0.,   0.5,   1. }
\definecolor{mg}{rgb}   {0.85,  0.,    0.85}
\definecolor{or}{rgb}   {0.9,  0.5,   0.}

\definecolor{webred}{rgb}{0.75,0,0}
\definecolor{webgreen}{rgb}{0,0.75,0}
\usepackage[citecolor=webgreen,colorlinks=true,linkcolor=webred]{hyperref}

\newtheorem{theorem}{Theorem}[section]

\newtheorem{remark}[theorem]{Remark}

\newtheorem{proposition}[theorem]{Proposition}

\newtheorem{example}{\it Example\/}

\newcommand{\ar}{\rightarrow}
\newcommand{\Ab}{{\bf A}}

\begin{document}

\title
{Inequalities for the lowest magnetic Neumann eigenvalue}

\author{S. Fournais}
\address[S. Fournais]{Aarhus University, Ny Munkegade 118, DK-8000 Aarhus~C, Denmark}
\email{fournais@math.au.dk}
\address[B. Helffer]{LMJL, Universit\'e de Nantes, 2 rue de la Houssini\`ere, 44322 Nantes Cedex France.}
\email{Bernard.Helffer@univ-nantes.fr}

\author{B. Helffer}
\date{\today}

\begin{abstract}
We study the ground state energy of the Neumann magnetic Laplacian on planar domains.
For a constant magnetic field we consider the question whether the disc maximizes this eigenvalue for fixed area. More generally, we discuss old and new bounds obtained on this problem.
\end{abstract}

\maketitle
\section{Introduction}

\subsection{The setup}
We consider an open set $\Omega \subset \mathbb R^2$ that is smooth, bounded and  connected. We denote by $A(\Omega)$ the area of $\Omega$, and define $R_{\Omega}$ to be the radius of the disc with the same area as $\Omega$, i.e.
\begin{align}
\pi R_\Omega ^2 = A(\Omega).
\end{align}
Let $\lambda_1^N(B,\Omega)$ be the ground state energy for the magnetic Neumann Laplacian on $\Omega$ with constant magnetic field of intensity $B\geq0$, i.e.
\begin{equation}\label{defHN}
H^N({\bf A} , B,\Omega):=
(-i\nabla + B {\mathbf A})^2\,
\end{equation}
 where $${\mathbf A} = \frac{1}{2} (-x_2, x_1)$$ (in particular $\nabla \times {\mathbf A} = 1$), and where we impose (magnetic) Neumann boundary conditions. 
 
Similarly, $\lambda_1^D(B,\Omega)$ will denote  the ground state energy in the case where we impose the Dirichlet boundary condition. We are interested in upper and lower bounds on these eigenvalues, universal or asymptotic in the two regimes $B \ar 0$ or $B \ar +\infty\,$. 
When considering lower bounds,  we  first mention the following  result obtained by L.~Erd\"os \cite{Er1} (in the spirit of the Faber-Krahn inequality for non-magnetic eigenvalues)
\begin{theorem}\label{thm:Erdos}
For any planar domain $\Omega$ and $B  >0$, we have:
\begin{equation}\label{ineq}
\lambda^D_1(B , \Omega) \geq \lambda^D_1 (B , D(0,R))\,,
\end{equation}
Moreover the equality in (\ref{ineq})  occurs if and  only if $\Omega = D (0,R_\Omega)$.
\end{theorem}

We would like to analyze a similar question for the Neumann magnetic Laplacian.
 $$
 \boxed{\mbox{ Question 1: For which $B>0$ do we have }  
\lambda^N_1(B,\Omega)  \leq   \lambda^N_1 (B , D(0,R_\Omega)) \quad ?}
$$
When $\Omega$ is assumed to be simply connected, our choice of ${\bf A}$ such that the magnetic field ${\rm curl\,} {\bf A} =1$ is not important because, by gauge invariance,  this spectral question depends only on the magnetic field.  We will discuss the non simply connected situation in Section \ref{ss5.3}.

\medskip

To analyze Question 1, we first look at the two asymptotic regimes $B \ar 0$ and $B \ar +\infty\,$.

\subsection{Weak magnetic field asymptotics}
By rather standard perturbation theory \cite[Proposition 1.5.2]{FH-b}, we have the following weak field asymptotics.
\begin{theorem}\label{Th1.2}
Let $\Omega \subset {\mathbb R}^2$ be smooth, bounded and  simply connected. 
There exists a constant $C_{\Omega}>0$ such that for all $B>0$
\begin{equation}\label{eq:1.2}
A(\Omega)^{-1} B^2 \int_{\Omega} |{\mathbf A}'|^2 \,dx - C_{\Omega}  B^4  \leq \lambda_1^N(B, \Omega) \leq A(\Omega)^{-1}  B^2 \int_{\Omega} |{\mathbf A}'|^2 \,dx\,,
\end{equation}
where the magnetic potential ${\bf A}'$ is the solution of 
\begin{equation}
 \nabla \times {\mathbf A}' = 1\,, \,\nabla \cdot {\mathbf A}' =0 \mbox{  and } {\mathbf A}' \cdot \nu =  0 \mbox{ on } \partial \Omega\,.
\end{equation}
\end{theorem}
Notice that in the case of the disc, we have ${\mathbf A}' = {\mathbf A}\,$.
 A weak version of  Question~1  above would consequently be:
 $$\boxed{
\mbox{ Question 2: Do we have } \int_{\Omega} |{\mathbf A}'|^2 \,dx \leq   \frac{1}{4} \int_{D(0,1/\sqrt{\pi})} r^2\,dx   \mbox{ if } A(\Omega)=1\, ?}
 $$
We will review the affirmative answer to Question 2 in Section~\ref{SmallB} below.

\subsection{Strong magnetic  field asymptotics}
For a smooth domain $\Omega$ and a point $P \in \partial \Omega$ we denote by $\kappa(P)$ the curvature of the boundary at $P$. We denote by $\kappa_{\rm max}(\Omega)$ the maximum value of $\kappa(P)$, $P \in \partial \Omega$.

In the limit where $B \rightarrow + \infty$, we have the following \cite[Theorem 8.3.2]{FH-b} (referring to former results by Bernoff-Sternberg \cite{BeSt}, Helffer-Morame \cite{HeMo}, Lu-Pan \cite{LP}).
\begin{theorem}\label{Th1.3}
Let $\Omega \subset {\mathbb R}^2$ be smooth and bounded. 
There exist $c_{\Omega}, B_{\Omega}>0$ such that
\begin{align}\label{expansion}
\left| \lambda_1^N(B, \Omega) - \left( \Theta_0 B - C_1\kappa_{\rm max}(\Omega) B^{1/2}\right)\right| \leq c_{\Omega} B^{1/3},
\end{align}
for all $B \geq B_{\Omega}$. Here $\Theta_0, C_1>0$ are universal constants, in particular, independent of $B$ and $\Omega$.
\end{theorem}

\begin{remark}
In the paper by Baumann-Phillips-Tang \cite{BaPhTa} (Theorem 6.1, p.~24)
 the authors prove the following more precise  asymptotic expansion for large values of $B$ in the case of the disc (see Fournais-Helffer  \cite[Chapter 5]{FH-b} and Fournais-Persson \cite{FPS} for improvements):
\begin{equation}\label{fPh}
 \lambda^{N} (B, D(0,R)) = \Theta_0 B  - C_1\, \frac 1 R \, B ^\frac 12 +
\mathcal O (1)\,.
\end{equation}
When $\kappa(P)$ has a unique non degenerate maximum, a more complete expansion than \eqref{expansion} can be obtained (see  Fournais-Helffer \cite{FH}).
\end{remark}

The asymptotics for strong magnetic fields leads us to the next question
$$\boxed{
\mbox{ Question 3: Is the maximal boundary curvature minimized by the disc }  ?}
 $$
We will review the affirmative answer to Question 3 (for simply connected domains) in Section~\ref{MaxCurv} below. 

\subsection{Reverse Faber-Krahn inequality for magnetic fields}
The analysis of Question 2 and 3, i.e. the study of the limits of large and small magnetic field strength,  
suggests that
\begin{align}\label{eq:FaberKrahnInverse}
\lambda^N_1(B,\Omega)  \leq   \lambda^N_1 (B , D(0,R_\Omega)),
\end{align}
for all $B$. This would correspond to a reverse Faber-Krahn inequality for magnetic fields, i.e. to an affirmative answer to Question 1.
Notice though, that we do not prove such an inequality in this paper. Also notice that this inequality is not true in general non-simply connected domains as the counterexample in Remark~\ref{rem:Counter} below shows. Notice also that it should not be confused with  the inequality---which often goes under the name `reverse Faber-Krahn inequality'
\begin{align}\label{SzegoWeinberger}
\lambda^N_2(0,\Omega)  \leq   \lambda^N_2 (0 , D(0,R_\Omega))\,,
\end{align}
which is due to Szeg\"o for the two-dimensional case
and Weinberger for the general case, see \cite{Sz, We}.

\begin{remark}
The discussion in the present paper also applies to a magnetic version of \eqref{SzegoWeinberger}. I.e. one may ask if the inequality
\begin{align}\label{SzegoWeinberger-magnetic}
\lambda^N_2(B,\Omega)  \leq   \lambda^N_2 (B , D(0,R_\Omega))\,
\end{align}
holds for all $B>0$. Since \eqref{SzegoWeinberger} is strict if $\Omega \neq D(0,R_\Omega)$ we immediately get \eqref{SzegoWeinberger-magnetic} for small values of $B$ by continuity. 
For $B$ large one can argue as follows: 
Suppose $\Omega$ is simply connected and not a disc.
By the affirmative answer to Question 3 given in Section~\ref{MaxCurv} below and continuity of the curvature we may choose distinct points $x, x' \in \partial \Omega$ such that the respective curvatures satisfy
$$
\kappa(x) \geq \kappa(x') > 1/R,
$$
(with $1/R$ being the curvature of the disc $D(0,R_\Omega)$). The proof of \eqref{expansion} (see \cite{HeMo}) involves the construction of approximate eigenfunctions localized near an arbitrary boundary point and with $\kappa_{\rm max}(\Omega)$ replaced by the boundary curvature at that point. Since the points $x, x'$ are distinct these approximate eigenfunctions will have disjoint support for large enough $B$ and both give energy expectations below the value \eqref{fPh} for the disc.
The inequality \eqref{SzegoWeinberger-magnetic} follows upon an application of the variational characterization of eigenvalues. Actually, by generalization to arbitrary number of points we find that if $\Omega$ is simply connected and not a disc, then for all $n \geq 1$ there exists $B_0>0$ such that
$$
\lambda^N_n(B,\Omega)  \leq   \lambda^N_1 (B , D(0,R_\Omega)), \qquad \text{ for all } B \geq B_0.
$$
However, to establish \eqref{SzegoWeinberger-magnetic} for intermediate values of $B$ remains open.
\end{remark}

\section{Around maximal curvature}\label{MaxCurv}
In this section we assume that $\Omega$ is simply connected.
We discuss the following result.
\begin{theorem}\label{thm:Pank}
For given area, the maximal curvature is minimized by the disc.
\end{theorem}
This is actually an old theorem which was rediscovered in  \cite{Pa} (and proved in the star-shaped case). 
\begin{proof}[Sketch of the proof of Theorem~\ref{thm:Pank} in the star-shaped case]
If $x(s)= (x_1(s), x_2(s))$ is a point of the boundary (parametrized by the arc length coordinate), it is an immediate consequence of the divergence theorem that
\begin{equation}\label{eq:area-p}
A(\Omega) = \frac 12 \int_{\partial \Omega} p(s) ds\,,
\end{equation}
where
$$
p(s)= x(s)\cdot \nu(s)\,,
$$
and $\nu(s)$ is the outward normal\footnote{Note that  in \cite{FH-b} the normal is directed inwards.}.
We note that if $\Omega$ is star-shaped with respect to $0$, then $p(s) \geq 0\,$.

\medskip

The second ingredient in the proof is the so-called Minkowski formula which in dimension $2$ reads
\begin{equation} \label{f2-new}
\ell (\partial \Omega) =  \int p(s) \kappa (s) ds\,.
\end{equation}

The identity \eqref{f2-new} follows from the observation that
$$
x''(s) = - \kappa (s) \nu(s)\,.
$$
Hence we can rewrite 
$$
\int p(s) \kappa (s) ds = - \int x(s) x''(s) ds = \int |x'(s)|^2 ds = \int ds=  \ell (\partial \Omega)\,.
$$
Combining \eqref{f2-new} with \eqref{eq:area-p} (and the positivity of $p(s)$), we easily get
\begin{equation}\label{bh1}
\ell (\partial \Omega)\leq  2 \kappa_{\rm max}(\Omega)   A(\Omega)\,.
\end{equation}
Now the classical isoperimetric inequality says that
\begin{equation}
\label{bh2}
4 \pi A(\Omega) \leq \ell(\partial \Omega)^2\,,
\end{equation}
with equality (only)  in the case of the disc.

Inserting \eqref{bh2} in \eqref{bh1} we find
\begin{equation}\label{eqpp}
\sqrt{\pi} A(\Omega)^{-\frac 12} \leq \kappa_{\rm max}(\Omega) \,.
\end{equation}
with strict inequality if $\Omega$ is not a disc.
\end{proof}

The general (simply connected but not star-shaped) case was open until recently. However, it has been settled in  \cite{Pa1} on the basis of a result due to Pestov-Ionin \cite{PI} (and \cite{HT}):
\begin{proposition}
For a smooth closed  Jordan curve, the interior of the curve contains a disk of radius $\frac{1}{\kappa_{\rm max} (\Omega)}$.
\end{proposition}

Finally, we mention that, as observed in  \cite{Pa}, this implies through the semi-classical analysis recalled in Theorem \ref{Th1.3}
 that in the large magnetic field strength limit we have
 $$
 \lambda_1^N(B, \Omega) \leq \lambda_1^N(B, D(0,R_\Omega))
  + {\mathcal O}(B^{1/3}),\qquad \text{ as } B \rightarrow \infty.
 $$
From this we deduce:
\begin{proposition} 
Let $\Omega \subset {\mathbb R}^2$ be smooth and simply connected. There exists $B_1(\Omega)  >0$ such that, for all $ B\geq B_1(\Omega)$, 
 \begin{equation}\label{eq:LargeFieldIneq}
 \lambda^N_1 (B, \Omega) \leq  \lambda^N_1 (B, D(0,R_\Omega))\,.
 \end{equation}
 Furthermore, the inequality \eqref{eq:LargeFieldIneq} is strict unless $\Omega = D(0,R_\Omega)$.
 \end{proposition}

\begin{remark}[Non-simply-connected counterexample to reverse Faber-Krahn]\label{rem:Counter}
Let $\Omega$ be an annulus, i.e. $\Omega = \{x\in{\mathbb R}^2\,: R_1 < |x| < R_2\}$ for some $0<R_1<R_2$. Notice that $\kappa_{\rm max}(\Omega) = \kappa_{\rm max}(D(0,R_2)) = \frac{1}{R_2}$, since the curvature on the inner boundary of the annulus is negative.
Therefore, $\kappa_{\rm max}(\Omega) < \kappa_{\rm max}(D(0,R_{\Omega}))$, and we get from Theorem~\ref{Th1.3} that
$$
 \lambda^N_1 (B, \Omega) >  \lambda^N_1 (B, D(0,R_\Omega)),
$$
for sufficiently large $B$.
\end{remark}

\section{Torsional rigidity}
\label{SmallB}
 We assume that $\Omega$ is simply connected and introduce 
\begin{equation}
\widehat S_\Omega:=\int_{\Omega} |{\mathbf A}'|^2 \,dx
\end{equation}
where the magnetic potential ${\bf A}'$ is the solution of 
\begin{equation}
 \nabla \times {\mathbf A}' = 1\,, \quad {\rm div\,}  {\mathbf A}' =0 \quad \mbox{  and } \quad {\mathbf A}' \cdot \nu =  0 \mbox{ on } \partial \Omega\,.
\end{equation}

As observed in \cite{FH-b} we have the identity
\begin{equation}
\widehat S_\Omega=\inf_{\phi \in H^2(\Omega)}  \int_\Omega  |\nabla \phi +{\bf A}'|^2 \,dx,
\end{equation}
which will be useful later.

Define now 
$\psi =\psi_\Omega$ to be the solution of
\begin{equation}\label{eq:3.7}
 \Delta \psi =1\,, \qquad \psi\big|_{\partial \Omega} =0\,.
 \end{equation} 
 Then we have 
 $${\bf A}' = \nabla^\perp \psi\,,
 $$
 where $ \nabla^\perp \psi = (-\partial_{x_2} \psi, \partial_{x_1} \psi)\,$.
 Hence, we get:
  $$
   \int_{\Omega} |{\mathbf A}'|^2 \,dx = \int_{\Omega}|\nabla \psi|^2 \,dx\,.
  $$
  The quantity 
  \begin{equation}\label{defsomega}
S_\Omega: =\int_{\Omega}|\nabla \psi|^2 \,dx
  \end{equation}
  with $\psi$ solution of \eqref{eq:3.7}, 
  is a well known quantity in Mechanics, which is called (see \cite{Sp}, up to a factor $2$) the torsional rigidity of $\Omega$.  In Mathematics, these quantities are analyzed in the celebrated book of Polya and Szeg\"o \cite{PS} where the results are obtained or illustrated  with many explicit computations for specific domains. By an integration by parts, we get 
  \begin{equation}\label{propomega}
S_\Omega = -  \int_\Omega \psi  \, dx\,.
\end{equation}

If $\psi =\psi_\Omega$ is the solution of  \eqref{eq:3.7}  then, by the maximum principle, $\psi < 0$ in $\Omega$ and attains its infimum $\psi_{min} (\Omega)$ in $\Omega$.

In \cite{HPS} it was observed, using Theorem~\ref{thm:Erdos} and  the asymptotics for $B$ large, that:
 \begin{equation}\label{converse}
0 >  \psi_{min} (\Omega) \geq \psi_{min} (D(0,R_\Omega))
 \end{equation}
 where $D(0,R)$ is the disk of same area as $\Omega$.
 As recalled in Sperb \cite[p. 193]{Sp} , this result   is actually due to Polya-Szeg\"o \cite{PS}. The inequality (3.7) is also a consequence of the very classical
result by G. Talenti, see \cite{Ta}.
   
To address Question 2 from the introduction, we will compare $S_{\Omega}$ for different domains.

\begin{example}[Optimizing $S_{\Omega}$ over ellipses]~\\ As a preliminary exercise, we can consider the case of the full 
 ellipse defined as $\{(x,y)\in \mathbb R^2\,,\, \frac{x^2}{a^2} + \frac{y^2}{b^2} < 1\}$, where explicit formulas are available. We have then 
 \begin{equation}\label{ellipse}
 \psi _{a,b}(x,y) = \frac{ 1}{ 2/a^2 + 2/b^2} \, \Bigl( \frac{x^2}{a^2} + \frac{y^2}{b^2} -1 \Bigr)\,,
 \qquad
 \psi_{a,b,min} = - \frac{ 1}{ 2/a^2 + 2/b^2} \,.
 \end{equation}
 To compare with equal area we assume in addition 
 $$
 ab =1\,,
 $$
 the unit disk corresponding to $a=b=1$.
 
 We then have by changing variables in the integral,
\begin{align*}
 - \int \psi_{a,b} (x,y) dx dy &= \frac{ 1}{ 2/a^2 + 2/b^2} \, \int \Bigl( 1- \frac{x^2}{a^2} - \frac{y^2}{b^2}  \Bigr)_+ dx dy \\
 & =  \frac{ 2ab }{ 1/a^2 + 1/b^2} ( - \int \psi_{1,1} (x,y) \,dx dy \,.
 \end{align*}
 This implies
 \begin{equation}
\left(  - \int \psi_{a,b} (x,y)\, dx dy\right) \leq \left( - \int \psi_{1,1} (x,y)\, dx dy\right)\qquad  \mbox{ if } ab =1\,,
\end{equation}
with equality for $a=b=1$.
\end{example}

This is actually the particular case (already mentioned in \cite{PS}) of  a general result communicated to us by D.~Bucur \cite{Bu} (see also \cite{BuGi} for generalizations to other models).
\begin{proposition}\label{prop:Bucur} Suppose that $\Omega$ is simply connected, then 
\begin{equation}
S_\Omega \leq S_{D(0,R)}\,.
\end{equation}
\end{proposition}

\begin{remark}
This problem was raised by Saint-Venant
as early as the 19th century (so this inequality is often called the Saint-Venant inequality) even though the proof is attributed
to G. Polya. 
\end{remark}

\begin{proof}[Sketch of the proof]
The proof is based on the formula
\begin{equation}
S_\Omega = - \left( \int_\Omega |\nabla \psi_\Omega|^2\, dx  + 2  \int_\Omega \psi_\Omega \, dx \right)
\end{equation}
which is an immediate consequence of \eqref{defsomega}. One can then follow the standard proof of the Faber-Krahn inequality.
Observe that $\psi_{\Omega}$ is the unique critical point
for the functional
$H_0^1(\Omega)\ni u \mapsto \int_\Omega |\nabla u|^2\, dx  + 2  \int_\Omega u \, dx$.
By uniqueness, it is a minimum. Using the Schwarz symmetrization procedure (see for example the survey \cite{As} for the main definitions and properties), this leads to 
$$  \int_\Omega |\nabla \psi_\Omega|^2\, dx  + 2  \int_\Omega \psi_\Omega \, dx \geq \int_{D(0,R)} |\nabla \psi^*_{\Omega}|^2\, dx  + 2  \int_{D(0,R)} \psi^{*}_\Omega \, dx
$$
where $\psi^*_\Omega$ is deduced from $\psi_\Omega$ by the symmetrization. One can then use that
 $$
 \int_{D(0,R)} |\nabla \psi^*_{\Omega}|^2\, dx  + 2  \int_{D(0,R)} \psi^{*}_\Omega \, dx \geq \int_{D(0,R)} |\nabla \psi_{D(0,R)}|^2\, dx  + 2  \int_{D(0,R)} \psi_{D(0,R)} \, dx\,.
 $$
\end{proof}

 As a corollary, we obtain
 \begin{proposition} 
 Suppose that $\Omega$ is smooth, bounded and simply connected.
 There exists $B_0=B_0(\Omega) >0$ such that, for all $ B \in (0,B_0)$, 
 \begin{equation}
 \lambda^N_1 (B, \Omega) \leq  \lambda^N_1 (B, D(0,R_\Omega))\,.
 \end{equation}
 \end{proposition}
 
 \begin{remark}
 Using recent results by Brasco-De Philippis-Velichkov \cite{BDV}, it is possible to show that we can take, assuming $ A(\Omega)=1$,
 $$ 
 B_0(\Omega) = C \mathcal A (\Omega)\,,
 $$
 where   $C>0$ is a universal constant and
 $ \mathcal A (\Omega)$ is the Fraenkel assymmetry 
 $$
 \mathcal A (\Omega):=\inf _{\text{ unit discs } D} A(\Omega \Delta D)\,,
 $$
 where the symbol $\Delta$ stands for the symmetric difference between sets.
 \end{remark}
 
 Observing that the main term of the asymptotics in Theorem~\ref{Th1.2} is an upper bound (coming from using a constant function as a trial state) we also get the following consequence of Prop.~\ref{prop:Bucur}.
 
 \begin{proposition} Suppose that $\Omega$ is simply connected. Then 
for all $ B >0$, 
 \begin{equation}\label{upperboundnsc}
 \lambda^N_1 (B, \Omega) \leq A(\Omega)^{-1}  S_\Omega B^2   \,.
 \end{equation}
 and
 \begin{equation}\label{est2}
 \lambda^N_1 (B, \Omega)  
   \leq \frac {A(\Omega)}{8\pi} B^2\,.
   \end{equation}
   
 \end{proposition}
 \begin{remark}
 Note also that \eqref{est2} can also be obtained via  Polya's inequality:
$$
A(\Omega)^{-1}  S_\Omega  \leq \frac{1}{\lambda_1^D(0,\Omega)}\,,
$$ 
(which results of a combination of \eqref{defsomega}, \eqref{propomega} and the Poincar\'e inequality)
and the standard Faber-Krahn inequality.

The recent improvement in \cite{BNFT} of the Polya inequality permits actually  to improve \eqref{est2}.
\end{remark}
\begin{remark}
One can find in \cite{PS} (see p.~10) a lower bound for smooth, simply connected planar domains:
\begin{equation}\label{eq:lb}
S_\Omega \geq  \frac{\pi}{8}  R_{in}^4\,,
\end{equation}
with optimality in the case of the disk.
\end{remark}

\section{lower bounds}
There are very few universal lower bounds in the literature.  Theorem 3.8 in the work of Ekholm-Kowarik-Portmann \cite{EKP}
 reads:
 \begin{theorem} Under the same assumptions as before
\begin{equation}
\label{eq:EKP1}
\lambda^N_1(B,\Omega) \geq \frac{\pi}{4 A(\Omega)} \, \frac{ B^2 R_{in}^4 \lambda_2^N(0,\Omega)}{B^2 R_{in}^2 + 6 \lambda_2^N (0,\Omega)} \,\mbox{ if } 0\leq B \leq R_{in}^{-2}\,,
\end{equation}
and
\begin{equation}\label{eq:EKP2}
\lambda^N_1(B,\Omega) \geq \frac{\pi}{32 A(\Omega)} \, \frac{ E(R_{in}^2 B )  \lambda_2^N(0,\Omega)}{B  + 12 \lambda_2^N (0,\Omega)} \,\mbox{ if } B \geq R_{in}^{-2}\,,
\end{equation}
where  $R_{in}$ is the interior radius of $\Omega$ and $E(t)$ is the largest integer $\leq t\,$.
\end{theorem}
Note that this estimate is coherent with the homogeneity property:
$$\lambda^N_1 (t^{-2}B, t\Omega) = t^{-2} \lambda^N_1(B,\Omega)\, , \, \forall t>0\,,$$
but none of these estimates are close to the associated asymptotics ($B$ small or $B$ large).

As $B \ar 0$, we indeed get from \eqref{eq:EKP1} that
\begin{equation}
\lambda^N_1(B,\Omega) \geq \frac{\pi}{24 A(\Omega)} \,  B^2 R_{in}^4 + O (B^3)\,,
\end{equation}
which should be compared to the asymptotic result using \eqref{eq:lb}. 

As $B\ar + \infty$, \eqref{eq:EKP2} becomes
\begin{equation}\label{eq:AssympLBound}
\lambda^N_1(B,\Omega) \geq \frac{\pi}{32 A(\Omega)} \, R_{in}^2   \lambda_2^N(0,\Omega) + o(1) \,,
\end{equation}
In particular, the right hand side of \eqref{eq:AssympLBound} remains bounded in contrast to \eqref{expansion}.
Other geometric  upper bounds in \cite{Er1} and \cite{CS,CEIS} seem to indicate that this lower bound is rather far from an optimal one for $B$ large.

It could be also interesting to compare with what one gets from Section 3.2 in our book \cite{FH}. 
\section{Extensions and open  questions. }
\subsection{Open questions} Of course, the main question is:
{\it Can we prove the reverse Faber-Krahn inequality \eqref{eq:FaberKrahnInverse} for arbitrary $B$ ?}

\medskip

Let us also mention the following connected questions
\begin{enumerate}
\item What can we say when $B$ is no more constant ?
\item What can we say in three dimensions ?
\item What can we say in the non-simply connected case ?
\end{enumerate}

In the next subsections we will discuss partial answers to these questions. 

\subsection{The case of a non constant magnetic field.}
Let us assume that the magnetic potential ${\bf A}$ has   as magnetic field $ \beta(x)$, with $\beta$ not necessarily constant.
Define,
$$
\widehat S_\Omega^{\bf A}:= \int_\Omega |{\bf A}'| ^2\, dx\,,
$$
where  $ {\bf A}'$ is the unique magnetic potential, such that $ {\bf A'} -{\bf A}$ is a gradient  and satisfying 
$$
{\rm curl\,} {\bf A}'=\beta,\quad  {\rm div\,} {\bf A}' =0 \mbox{ in } \Omega \qquad \mbox{ and }  \qquad {\bf A}'\cdot \nu =0  \mbox{ on } \partial \Omega\,.
$$
We have the following easy perturbation proposition:
\begin{proposition}\label{propnew} If $\Omega$ is connected, 
\begin{equation}
A(\Omega)^{-1} B^2 \widehat S_\Omega^{\bf A} - C_\Omega B^4 \leq \lambda_1^N (B {\bf A} , \Omega) \leq A(\Omega)^{-1} B^2 \widehat S_\Omega^{\bf A} \,.
\end{equation}
\end{proposition}
We now asssume that $\Omega$ is simply connected  (we recall that when $\Omega$ is simply connected, $\widehat S_\Omega^{\bf A}$ depends only on $\beta$) and 
define,
$$
S_\Omega^\beta:= \int_\Omega |\nabla \psi| ^2\, dx\,,
$$
where this time $\psi (x) $ is the solution of
$$
\Delta \psi (x) = \beta (x) \quad \mbox{ in } \Omega\,,\qquad \psi  =0 \quad \mbox{ on } \partial \Omega\,.
$$
We observe as in \cite{HPS1}, that ${\bf A}'=\nabla^\perp \psi\,$, which implies:
\begin{equation}
\widehat S_\Omega^{\bf A} = S_\Omega^\beta\,.
\end{equation}

We would like now to find an isoperimetric inequality for $S_\Omega^\beta$. From now on, we also assume  that the magnetic field $\beta$ satisfies
$$
\beta(x) >0\quad\mbox{ on } \overline{\Omega}\,.
$$

We then conclude from the maximum principle that $\psi < 0$ in $\Omega$ and follow the different steps of the constant magnetic field case.

First we observe that by an integration by parts, we can rewrite $S_\Omega^\beta$ in the following form, if $\Omega$ is simply-connected, 
\begin{equation}\label{eq:5.2}
S_\Omega^\beta = - \int_\Omega  \beta(x) \psi (x)\,.
\end{equation}
For applying the variational argument of D. Bucur, we rewrite $S_\Omega^\beta$ in the form
$$
S_\Omega^\beta  = - \left (  \int_\Omega |\nabla \psi (x) |^2\, dx  + 2  \int_\Omega \beta (x) \psi (x)  \, dx    \right)\,.
$$
We now observe that by the Schwarz symmetrization procedure, we get (see \cite[p.~4, lines -2 and -1]{As}).
$$ 
\int_\Omega |\nabla \psi (x) |^2\, dx\,  +\,  2  \int_\Omega \beta (x) \psi (x)  \, dx  
\geq  \int_{D(0,R_\Omega)} |\nabla \psi^* (x) |^2\, dx\,  + \, 2  \int_{D(0,R_\Omega)} \beta^* (x) \psi^* (x)  \, dx \,,
$$ 
with $\psi^*\in H_0^1(\Omega)$.

Here we have used \cite[Theorem 3.4]{LL} (or \cite[(1.14)]{As}). Hence we get
\begin{proposition}
If $\Omega$ is simply connected and $\beta\geq 0$ then 
\begin{equation}
S_\Omega^\beta \leq S_{D(0,R_\Omega)}^{\beta^*}\,,
\end{equation}
where $\beta^*$ is the Schwarz symmetrization of $\beta$.
\end{proposition}

\begin{remark} This proposition is also a consequence of Talenti's result,  see  \cite{Ta}.
\end{remark}

We can  then be a little more explicit. Since $\beta^*$ is radial\footnote{In the case of the radial function, we use the same notation for the function and the corresponding $1D$ function.} it is possible to compute the corresponding $ S_{D(0,R_\Omega)}^{\beta^*}$.
We have first to compute the (radial) solution  $ \psi^*$ of
$$
\Delta \psi (x) = \beta ^* (x) \mbox{ in } D(0,R_\Omega)\,,\, \psi =0 \mbox{ on } \partial D(0,R_ \Omega)\,.
$$
 Hence, we have to analyze the solutions in the interval $(0,R)$  of
\begin{equation} 
\frac{d^2\psi }{dr^2} + \frac 1r \frac{d\psi }{dr} =  \beta^*(r)\,, \, \psi(R)=0\,.
\end{equation}
This is rewritten in the form
\begin{equation}
 \frac 1r  (r\psi')' (r) = \beta^*(r) \mbox{ in } (0,R)\,, \, \psi (R)=0\,.
\end{equation}
We get first
$$
r \psi' (r) =  \int_ 0  ^r  t \, \beta^*(t)\, dt + C\,,
$$
and then (we add the condition that $\psi$ is continuous at $0$) we get 
\begin{equation}
\psi^* (r) = - \int_r^R \frac{1}{s}   (\int_0 ^s  t \, \beta^*(t)\, dt ) \,ds
\end{equation}
Hence we finally have
$$
S_{D(0,R_\Omega)}^{\beta^*} = 2 \pi \int_0^R r \beta^*(r) \left( \int_r^R \frac{1}{s}   \left(\int_0 ^s  t \, \beta^*(t)\, dt \right) \,ds\right) dr\,.
$$
 Define the flux function
$$
F(r) = \int_0^r r' \beta^*(r')\,dr'.
$$
By an integration by parts, we then find,
$$
S_{D(0,R_\Omega)}^{\beta^*} = 2 \pi \int_0^{R_\Omega} \frac{F(r)^2}{r} \,dr.
$$

\begin{remark}\label{rem:5.3}
In our particular case, we recover from Proposition \ref{propnew},  under the additional assumption that $\Omega$ is simply connected, the following variant\footnote{ The authors use another lowest eigenvalue corresponding to a Laplacian on $2$-forms  satisfying specific boundary conditions   but work in any dimension. Here we are in dimension 2 and identify $2$-forms and functions.} of Colbois-Savo's result (Proposition 4 in \cite{CS} or \cite{CEIS}):
\begin{equation}
\lambda_1^N (B{\bf A}, \Omega)  \leq  B^2 \frac{ 1}{A(\Omega)\lambda_1^D(\Omega)} \left(\int_\Omega \beta^2 dx \right)\,,
\end{equation}
where $\lambda_1^{D}(\Omega)=\lambda_1^{D}(0,\Omega)$ is  the groundstate energy of the Dirichlet realization of the Laplacian in $\Omega$.

The authors use there (with $k=1$)  the comparison theorem between the magnetic  Laplacian and a Schr\"odinger operator with electric potential: $V(x)=B^2 A(x)^2$
\begin{equation}
\lambda_k^N (B{\bf A}, \Omega)\leq \lambda_k^N (-\Delta + V, \Omega)\,,\, \forall k \geq 1\,,
\end{equation}
which results from the min-max comparison principle and is true under the assumption that $\Omega$ is connected.\\
Here we observe, using \eqref{eq:5.2}, that, if $\Omega$ is simply connected, 
$$
S_\Omega^\beta  \leq \|\beta\| \, \|\psi\| \leq \frac{ \|\beta\| }{\sqrt{\lambda_1^D (\Omega)}} \left(S_\Omega^\beta\right)^\frac 12\,,
$$
which gives, by the standard isoperimetric inequality for $\lambda_1^D$
\begin{equation}
S_\Omega^\beta  \leq \|\beta\| ^2 \lambda_1^D(\Omega)^{-1} \leq \|\beta\| ^2 \lambda_1^D(D(0,R_\Omega))^{-1}\,.
\end{equation}
\end{remark}
\begin{remark}
As in the constant magnetic field case, the previous estimates are only good for small values of $B$. When $B$ is large, we refer to the semi-classical analysis of N. Raymond \cite{Ra} or to the universal estimates of \cite{Er1} or \cite{CS,CEIS}.
\end{remark}

\subsection{The 3D-case.}~\\
 We consider the ground state energy $\lambda^N (\Omega,B) $ of the Neumann magnetic Laplacian in $\Omega$ attached to the magnetic potential $B(-\frac{x_2}{2}, \frac{x_1}{2}, 0)$. 
There is no hope to have in 3D the inequality
\begin{equation}\label{ineq3D}
\lambda^N (\Omega,B)  \leq \lambda^N (B (0,R_\Omega), B ) \quad \text{ if }  \quad |\Omega| = |B(0, R_\Omega)| = \frac 43 \pi R_\Omega^3\,,
\end{equation}
where $B(0,R)$ denotes the ball of radius $R$ centered at $0$ in $\mathbb R^3$  and $|\Omega|$ denotes the volume of $\Omega\,$.\\
Take indeed, for some $L>0$, $\omega \subset {\mathbb R}^2$, the set $\Omega_L= \omega \times [0,L]\,$. We have 
$$
|\Omega_L| = A(\omega)\, L =\frac{4\pi R_{\Omega_L}^3}{3} \,,
$$
and (by separation of variables)
$$
\lambda^N(\Omega^L, B) = \lambda^N (\omega,B)\,.
$$
But, using the constant function $|\Omega_L|^{-\frac 12}$ as trial state, we have, for some universal constant $C>0\,$, 
$$
\lambda^N (B (0,R_{\Omega_L}), B) \leq C B^2 R_{\Omega_L}^2\,.
$$
Taking the limit as $L \rightarrow 0$, we see that \eqref{ineq3D} cannot be satisfied.

This is actually not surprising because the magnetic field introduces a privileged direction. The "optimal domain" should have the same property.
\subsection{The nonsimply connected case}\label{ss5.3}~\\
We mention a recent preprint of B.~Colbois and A.~Savo \cite{CS}, further developed\footnote{We thank B. Colbois for communicating to us these papers before publication.} in \cite{CS1} for the lower bounds and \cite{CEIS} for the upper bounds,  devoted to the Neumann problem (see in particular the upper bound given in their Proposition 4  in \cite{CS} and a lower bound when $B=0$) and two papers by Helffer-Persson Sundqvist \cite{HPS,HPS1} initially motivated by Ekholm-Kowarik-Portman \cite{EKP}. Here we denote by $\lambda^N ({\bf A}, \Omega)$  the first eigenvalue of the Neumann problem. We are not aware 
 of a comparison with a specific domain in the non simply connected case.
 
 As in \cite{HPS1} (see also \cite{CS} for a more geometric formalism), we can observe that  if $\Omega$ has $k$ holes $D_j$ ($j=1,\cdots,k$), 
 $$
 \Omega:= \widetilde \Omega \setminus \cup_j D_j\,,
 $$
 and $\Ab$ is a solution of
 $$
{\rm curl\,} {\bf A}=\beta, \quad {\rm div\,} {\bf A} =0 \quad \text{ in } \Omega, \qquad \text{ and }  \qquad {\bf A}\cdot \nu =0  \quad \text{ on } \partial \Omega\,,
$$
  the generating function  $\psi$ such that ${\bf A}=\nabla^{\perp}\psi$ 
  is now solution of 
  \begin{equation}
  \Delta \psi=\beta \quad \text{ in } \quad \Omega, \qquad \text{ with } \psi |_{\partial \widetilde \Omega}=0  \quad \text{ and } \quad \psi |_{\partial D_j} = C_j\,,
  \end{equation}
  for some real constants $C_j$.
  
 \medskip
  
 We can then write 
 $$
 \psi =\psi^0 + \sum_{j=1}^k\, C_j \, \theta^j
 $$
 where $\psi^0$ is the solution of 
 $$
 \Delta \psi^0 = \beta, \qquad  \psi^0 |_{\partial \Omega} =0\,,
 $$
 and $\theta^j$'s is the solution of
 $$
 \Delta \theta^j =0, \qquad  \theta^j |_{\partial D_i} = \delta_{ij}\,,\, \theta^j |_{\partial \widetilde \Omega} =0\,.
 $$
   It remains to compute $S_\Omega^\beta$, which gives
 \begin{equation}\label{f1}
 S_\Omega^\beta = - \int_\Omega  \beta \psi^0 +  \sum_{i,j} C_iC_j M_{ij}  \,,
 \end{equation}
 with   
 \begin{equation}\label{f2}
  M_{ij}=  \int_\Omega  \nabla \theta^i\cdot \nabla \theta^j dx = M_{ji} \,.
  \end{equation}

We obtain by a reasoning similar to Remark~\ref{rem:5.3},
 $$
  S_\Omega^\beta \leq  \|\beta\|^2 \, \lambda_D(\Omega)^{-1}   + \sum_{i,j} C_iC_j  M_{ij} \,.
  $$
 When coming back to the upper bound to $\lambda_1^N (\beta B, \Omega)$, we have to implement the gauge invariance in order to minimize the $C_j$ as in \cite{HPS1}. 
  
 Here, it is better to return to the formulation in terms of the circulations of $\bf A$ along the $D_i$ ($i=1,\cdots, k$). We introduce
  $$
  {\bf A}^0 =\nabla^\perp \psi^0\,$$
  and the circulations of ${\bf A}^0$ and ${\bf A}$ along $\partial D_i$
  $$
  \Phi^0_i = \int_{\partial D_i}\, A_0 \, ds, \qquad \Phi_i  =  \int_{\partial D_i}\, {\bf A} \, ds \quad \text{ for } \text i=1,\ldots, k\,.
  $$
  
  Observing that
    $$
  M_{ij}= - \int_{\partial D_i} \partial_\nu \theta^j\,,
  $$
we get 
  \begin{equation}
  \Phi_i = \Phi^0_i - \sum_j \, M_{ij} C_j\,.
  \end{equation}
  We also note  that $M$ is positive definite 
  (just compute $M \vec C \cdot \vec C$ to get the injectivity). We deduce from \eqref{f1}:
  $$
 S_\Omega^\beta = - \int_\Omega  \beta \psi^0\, dx  +  \,| M^{-\frac 12} ( \Phi-\Phi^0) | ^2\,.
 $$
 Note here that, for fixed $\beta \geq 0 $, $\psi^0 \leq 0$ by the maximum principle and  $$- \int_\Omega  \beta \psi^0\, dx   \geq 0\,.$$ Hence for fixed $\beta \geq 0 $ the torsion rigidity is minimal when $\Phi=\Phi^0$.

\medskip

 Coming back to the upper bound for $\lambda^N ( {\bf A},\Omega)$, we can use the gauge invariance of the problem (see for example Proposition 2.1.3 in \cite{FH}) in order to get:
 \begin{equation}
 \lambda^N ( {\bf A},\Omega) \leq \frac{1}{A(\Omega)} \left( \|\beta\|^2 \lambda_1^D(\Omega)^{-1} + \inf_{\gamma \in \mathbb Z^k}  \left( \,|M^{-\frac  12} (\Phi-\Phi^0 - 2 \pi \gamma)| ^2\right)\right)
 \end{equation}
  
 \begin{remark}
 We can use the isoperimetric inequality for $\lambda_1^D$ and get
  \begin{equation}
 \lambda^N ( {\bf A},\Omega) \leq \frac{1}{A(\Omega)} \left( \|\beta\|^2 \lambda_1^D(D(0,R_\Omega))^{-1} + \inf_{\gamma \in \mathbb Z^k}  \left( \,|M^{-\frac 12} (\Phi-\Phi^0 - 2 \pi \gamma)| ^2\right)\right)
 \end{equation}
 \end{remark}

\begin{remark}
 These estimates are related to the so-called Aharonov-Bohm effect (see  references in \cite{CEIS}).
 \end{remark}

~\\
{\bf Acknowledgements}
The discussion on this problem started at a nice meeting in Oberwolfach (December 2014) organized by V. Bonnaillie-No\"el, H. Kova\v{r}\'ik and K.~Pankrashkin.

We would like to thank M. van den Berg,  D. Bucur, B. Colbois, M. Persson Sundqvist, K. Pankrashkin,  N.~Popoff and A. Savo  for discussions in the last two years around this problem. We also thank an anonymus expert for additional references.


\begin{thebibliography}{99}


\bibitem{As} M.S. Ashbaugh.
\newblock Isoperimetric and universal inequalities for eigenvalues.
\newblock London Math. Soc. Lecture Note Ser. 273, 95--139 (2000).

\bibitem{Ba} K. Ball.
\newblock Volume ratios and a reverse isoperimetric inequality. 
\newblock J. London Math. Soc. (2) 44 (1991), no. 2, 351--359.

\bibitem {BaPhTa} P.~Bauman, D.~Phillips, and Q.~Tang.
\newblock Stable nucleation for the Ginzburg-Landau system with an applied magnetic
 field.
\newblock Arch. Rational Mech. Anal. 142 (1998), 1--43.

\bibitem{BNFT} M. van den Berg, V. Ferone, C. Nitsch and C. Trombetti.
\newblock On Polya's Inequality for torsional rigidity
and first Dirichlet eigenvalue.
\newblock Integr. Equ. Oper. Theory 86 (2016), 579--600.


\bibitem {BeSt} A.~Bernoff and P.~Sternberg.
\newblock Onset of superconductivity in decreasing 
fields for general domains.
\newblock J.~Math.~Phys. 39 (1998), 1272--1284.

 \bibitem{BDV} L. Brasco, G. De Philippis, and B. Velichkov.
\newblock Faber-Krahn inequality in sharp  quantitative form.
\newblock arXiv:1306.0392v1 (June 2013).


\bibitem{Bu} D. Bucur.
\newblock Personal communication (March 2017).

\bibitem{BuGi} D. Bucur and A. Giacomini.
\newblock Faber-Krahn inequalities for the Robin-Laplacian: a free discontinuity approach. 
\newblock  Arch. Ration. Mech. Anal. 218, 757-824 (2015). 

\bibitem{CEIS} B. Colbois, A. El Soufi, S. Ilias, and A. Savo.
\newblock Eigenvalues upper bounds for the magnetic operator.
\newblock ArXiv:1709.09482v1, 27 Sep 2017.

\bibitem{CS} B. Colbois and A. Savo.
\newblock Eigenvalue bounds for the magnetic Laplacian.
\newblock ArXiv:1611.01930v1 (2016).

\bibitem{CS1} B. Colbois and A. Savo.
\newblock Lower bounds for the first eigenvalue of the magnetic Laplacian.
\newblock ArXiv:1709.09506v1 [math.DG] 27 Sep 2017.

\bibitem{EKP} 
T. Ekholm, H. Kova\v{r}\'ik, and F. Portmann.
\newblock Estimates for the lowest eigenvalue of magnetic Laplacians.
 \newblock  J. Math. Anal. Appl. 439 (1) (2016),  330--346. 

\bibitem {Er1} L.~Erd\"os.
\newblock Rayleigh-type isoperimetric inequality with a homogeneous magnetic field.
\newblock Calc. Var. and PDE. 4 (1996) , 283-292.

\bibitem{Er2} L. Erd\"os.
\newblock Spectral shift and multiplicity of the first eigenvalue of the magnetic Schr\"odinger operator in two dimensions.
\newblock Ann. Inst. Fourier, 52:1833-1874 (2002).

\bibitem{FNT}  V.~Ferone, C. Nitsch, and C.~Trombetti.
\newblock  On the maximal
mean curvature of a smooth surface. 
\newblock C. R. Math. Acad. Sci. Paris 354 (2016), no. 9, 891--895.


\bibitem{FH} S. Fournais, B. Helffer. 
\newblock Accurate eigenvalue asymptotics
for the magnetic Neumann  Laplacian. 
\newblock  {\it Ann. Inst. Fourier.} {\bf
  56} (1) 1-67 (2006).


\bibitem{FH-b}  S. Fournais, B. Helffer. 
\newblock  {\it Spectral Methods in
Surface Superconductivity.} 
\newblock Progress in Nonlinear Differential
Equations and Their Applications, Vol.~77, Birkh\"auser (2010).
  
  \bibitem{FPS} S. Fournais and M. Persson Sundqvist.
\newblock Lack of diamagnetism and the Little-Parks effect.
 \newblock  Comm. Math. Phys. 337 (1), 2015, 191--224. 
 
 \bibitem{HeMo} B. Helffer and A. Morame.
 \newblock Magnetic bottles in connection with superconductivity. 
 \newblock  Journal of Functional Analysis Vol. 185, No 2, October,  604--680 (2001). 
 

 \bibitem{HPS} B. Helffer and M. Persson Sundqvist.
 \newblock On the semi-classical analysis of the Dirichlet Pauli operator.
 \newblock J. Math. Anal. Appl. 449 (1), 2017, 138--153.

\bibitem{HPS1} B. Helffer and M. Persson Sundqvist.
 \newblock On the semi-classical analysis of the Dirichlet Pauli operator--  the non simply connected case.
 \newblock Problems in Mathematical Analysis and Journal of Mathematical Sciences, Vol. 226, No. 4, October, 2017. 
 
 \bibitem{HT} R. Howard and A. Treibergs.
 \newblock A reverse isoperimetric inequality, stability and extremal theorems for plane curves with bounded curvature.
 \newblock Rocky Mountain J. Math. 25 (1995) 635--684.
 
 \bibitem{Kaw} B. Kawohl. Overdetermined problems and the p-Laplacian. Acta
Math. Univ. Comenianae, 76 (2007), 77--83.
 
 \bibitem{LL} E.H. Lieb and M. Loss.
 \newblock Analysis.
 \newblock Graduate Studies in Mathematics 14. American Mathematical Society.
 
 \bibitem{LP} K. Lu and X. Pan.
 \newblock Eigenvalue problems of Ginzburg-Landau operator in bounded domains.
 \newblock J. Math. Phys. 40 (6) (1999) 2647--2670.

\bibitem{Pa} K. Pankrashkin.
\newblock An inequality for the maximum curvature of planar curves with applications to some eigenvalue problem.
\newblock https://arxiv.org/abs/1501.03792v3.

\bibitem{Pa1} K. Pankrashkin.
\newblock An inequality for the maximum curvature through a geometric flow.
\newblock  Arch. Math. 105 (2015), 297--300 (Springer).

\bibitem{PI} G. Pestov and V. Ionin.
\newblock On the largest possible circle embedded in a given closed curve.
\newblock Dokl. Akad. Nauk SSSR 127 (1959) 1170-1172 (in russian).

\bibitem{PS} G. Polya and G. Szeg\"o.
\newblock Isoperimetric Inequalities in Mathematical Physics.
\newblock Princeton University Press, Princeton, New Jersey (1951).

\bibitem{Ra} N. Raymond.
\newblock Sharp  asymptotics  for  the  Neumann  Laplacian  with  variable  magnetic   field  in  dimension 2.
\newblock Annales Henri Poincar\'e, 10(1), 95-122, (2009).


\bibitem{Sp} R. Sperb.
\newblock  Maximum principles and their applications.
\newblock Academic Press, New York, 1981.

\bibitem{Sz} G. Szeg\"o.
\newblock  Inequalities for certain eigenvalues of a membrane of given
area. 
\newblock J. Rational Mech. Anal. 3, (1954). 343-356.

\bibitem{Ta} G. Talenti.
\newblock Elliptic equations and rearrangements. 
\newblock Ann. Scuola
Norm. Sup. Pisa Cl. Sci. (4) 3 (1976), no. 4, 697--718.

\bibitem{We} H.F. Weinberger.
\newblock An isoperimetric inequality for the N-dimensional
free membrane problem. 
\newblock J. Rational Mech. Anal. 5 (1956), 633--636.



\end{thebibliography}
\end{document}